\documentclass[letterpaper]{article}

\usepackage{amsmath}
\usepackage{amssymb}
\usepackage{amsthm}
\usepackage{thmtools}
\usepackage{float}
\usepackage{mathtools}
\usepackage{enumitem}
\usepackage{booktabs}
\usepackage{microtype}

\usepackage[backend=biber, style=numeric, sorting=nyt, urldate=iso]{biblatex}
\AtEveryBibitem{
    \clearfield{note}
    \clearfield{eprint}
}
\DeclareSourcemap{
  \maps[datatype=bibtex]{
    \map{
      \step[fieldsource=doi,final]
      \step[fieldset=url,null]
      \step[fieldset=isbn,null]
    }
  }
}
\usepackage[hypertexnames=false]{hyperref}
\hypersetup{pdftitle={Two Perfect One-Factorizations of K92}, pdfauthor={Tomasz Maciosowski}}
\usepackage[capitalise,noabbrev]{cleveref}
\usepackage{orcidlink}

\declaretheorem[name=Theorem, numberwithin=section]{theorem}
\declaretheorem[name=Lemma, sharenumber=theorem]{lemma}

\declaretheorem[name=Example, style=definition, sharenumber=theorem]{example}

\crefname{theorem}{Theorem}{Theorems}
\crefname{lemma}{Lemma}{Lemmas}
\crefname{definition}{Definition}{Definitions}
\crefname{example}{Example}{Examples}

\title{Two Perfect One-Factorizations of $K_{92}$}
\author{Tomasz Maciosowski\,\orcidlink{0009-0000-7485-0786}\\
  \small Department of Mathematics and Statistics\\
  \small Memorial University of Newfoundland\\
  \small Canada
}

\date{}

\begin{document}
\maketitle

\begin{abstract}
  \noindent%
  A perfect 1-factorization of the complete graph $K_{2t}$ is a partition of its edges into perfect matchings such that the union of any two of them is a Hamilton cycle.
  Kotzig conjectured in 1963 that $K_{2t}$ has a perfect 1-factorization for every $t \ge 2$, but apart from some infinite families and a list of sporadic orders the conjecture is open, and $92$ was one of the nine orders up to $100$ for which no perfect 1-factorization was known.
  In this paper we construct two non-isomorphic perfect 1-factorizations of $K_{92}$.
\end{abstract}

\section{Introduction}\label{sec:intro}

A \emph{perfect 1-factorization} of the complete graph $K_{2t}$ is a partition of its edges into perfect matchings such that any two of them together form a Hamilton cycle.
Kotzig~\cite[p.~80]{kotzig:HamiltonGraphsHamilton:1964} conjectured in 1963 that $K_{2t}$ has a perfect 1-factorization for every $t \ge 2$; see the survey of Rosa~\cite{rosa:Perfect1factorizations:2019} for the history of the problem.
Beyond some infinite families, perfect 1-factorizations are known only for a list of sporadic orders, most of them of the form $q + 1$ with $q$ a prime power and found by computer searches for starters in the finite field $\mathbb{F}_q$ with prescribed multiplicative symmetries~\cite{dinitz:NewPerfectOnefactorizations:1989}; see~\cite{rosa:Perfect1factorizations:2019} and the list maintained by Wanless~\cite{wanless:Perfect1factorisations:}.
The smallest orders were settled one at a time, most recently $52$ by Wolfe~\cite{wolfe:PerfectOnefactorizationK52:2009} and $56$ by Pike~\cite{pike:PerfectOnefactorisation:2019}, and the orders up to $100$ for which the conjecture is open are $64$, $66$, $70$, $76$, $78$, $88$, $92$, $96$ and $100$.

In this paper we settle the case $2t = 92$ by constructing two non-isomorphic perfect 1-factorizations of $K_{92}$ (\cref{thm:main}).
Both are induced by starters in $\mathbb{Z}_{91}$ that are invariant under multiplication by a subgroup $U$ of order $3$ of the group of units $\mathbb{Z}_{91}^{*}$.
This is the idea behind the quotient coset starters of Dinitz and Stinson~\cite{dinitz:NewPerfectOnefactorizations:1989}, who used starters in a finite field $\mathbb{F}_q$ that are invariant under a subgroup of $\mathbb{F}_q^{*}$.
Since $91$ is not a prime power, there is no field with $91$ elements, so we work in the ring $\mathbb{Z}_{91}$ instead.
The same method was used by Ihrig, Seah and Stinson~\cite{ihrig:PerfectOnefactorizationK50:1987}, who found the first perfect 1-factorizations of $K_{50}$ among the starters in the ring $\mathbb{Z}_{49}$ that are invariant under multiplication by the subgroup $\{1, 18, 30\}$ of $\mathbb{Z}_{49}^{*}$.

The symmetry helps in two ways.
It cuts the number of pairs of factors whose union has to be checked from $45$ to $15$, and it turns the search for a starter into an exact cover problem.
It also restricts the choice of $U$: only subgroups whose non-identity elements fix no nonzero element of $\mathbb{Z}_{2t - 1}$ can lead to a perfect 1-factorization, and this limits the orders to which the method applies.

\section{Main Result}\label{sec:main}

\begin{theorem}\label{thm:main}
  The graph $K_{92}$ has two non-isomorphic perfect 1-factorizations.
\end{theorem}

The rest of this section develops the tools for the proof, which is given at the end of \cref{sec:91}.

\subsection{Multiplier Symmetry}\label{sec:multiplier}

A \emph{starter} in $\mathbb{Z}_n$ for an odd $n \ge 3$ is a set $S$ of $\frac{n - 1}{2}$ disjoint pairs $\{x, y\} \subset \mathbb{Z}_n \setminus \{0\}$ whose differences $\pm(x - y)$ cover $\mathbb{Z}_n \setminus \{0\}$ exactly once.
Starters are a standard tool for the construction of 1-factorizations~\cite[p.~622]{colbourn:HandbookCombinatorialDesigns:2006}.
A starter induces a 1-factorization of $K_{n + 1}$ on the vertex set $V = \mathbb{Z}_n \cup \{\infty\}$ with factors
\[F_i = (S + i) \cup \{\{i, \infty \}\} \quad \text{for } i \in \mathbb{Z}_n.\]
Since the translation $x \mapsto x + 1$, $\infty \mapsto \infty$ maps $F_i$ to $F_{i + 1}$, we have $F_i \cup F_j \cong F_0 \cup F_{j - i}$, and the translation $x \mapsto x - d$ gives $F_0 \cup F_d \cong F_0 \cup F_{-d}$.
Thus the 1-factorization is perfect if and only if $F_0 \cup F_d$ is a Hamilton cycle for all $d \in \{1, \ldots, \frac{n - 1}{2}\}$.
In the $K_{92}$ case this yields 45 conditions.
To reduce both the search space and the number of conditions, we look only for starters that are invariant under multiplication by the elements of a subgroup of the group of units $\mathbb{Z}_n^{*}$ of $\mathbb{Z}_n$.
Let $U$ be a subgroup of $\mathbb{Z}_n^{*}$.
Then $U$ acts on $\mathbb{Z}_n$ by $x \mapsto ux$ and on pairs by $u\{x, y\} = \{ux, uy\}$.
The \emph{$U$-orbit} of $x \in \mathbb{Z}_n$ is $Ux = \{ux \mid u \in U\}$ and the $U$-orbit of a pair $\{x, y\}$ is $\{\{ux, uy\} \mid u \in U\}$.
A set of pairs $S$ is \emph{$U$-invariant} if $uS = S$ for all $u \in U$, where $uS = \{u\{x, y\} \mid \{x, y\} \in S\}$.
Since the orbits of pairs partition the set of all pairs, $S$ is $U$-invariant if and only if it is a union of $U$-orbits of pairs.
Invariance says that $U$ fixes $S$ as a set, not that it fixes the pairs of $S$: each $u \in U$ permutes the pairs of $S$ among themselves.

\begin{example}\label{ex:orbits-91}
  Let $n = 91$ and let $U = \langle 16 \rangle = \{1, 16, 74\}$ be the subgroup generated by $16$, where $16^2 = 256 \equiv 74$ and $16 \cdot 74 = 1184 \equiv 1 \pmod{91}$.
  The $U$-orbit of $1$ is $\{1, 16, 74\}$ and the $U$-orbit of the pair $\{1, 41\}$ is
  \[
    \{1, 41\},\quad \{16, 16 \cdot 41\} = \{16, 19\},\quad \{74, 74 \cdot 41\} = \{74, 31\}.
  \]
  So a $U$-invariant starter that contains $\{1, 41\}$ also contains $\{16, 19\}$ and $\{31, 74\}$.
\end{example}

For $u \in \mathbb{Z}_n^{*}$ we define $\varphi_u \colon V \to V$ by $\varphi_u(x) = ux$ for $x \in \mathbb{Z}_n$ and $\varphi_u(\infty) = \infty$.
\Cref{lem:multiplier-automorphism,lem:orbit-conditions} are the analogues for $\mathbb{Z}_n$ of~\cite[Theorem~2.2 on p.~407 and Corollary~2.3 on p.~408]{dinitz:NewPerfectOnefactorizations:1989}, which concern starters in $\mathbb{F}_q$ that are invariant under the largest subgroup of $\mathbb{F}_q^{*}$ of odd order.

\begin{lemma}\label{lem:multiplier-automorphism}
  If $S$ is a starter in $\mathbb{Z}_n$ with induced factors $F_i$ and $u \in \mathbb{Z}_n^{*}$ satisfies $uS = S$, then $\varphi_u$ is an automorphism of $K_{n + 1}$ with $\varphi_u(F_i) = F_{ui}$ for all $i \in \mathbb{Z}_n$, and for every $d \in \mathbb{Z}_n$ the union $F_0 \cup F_d$ is a Hamilton cycle if and only if $F_0 \cup F_{ud}$ is.
\end{lemma}
\begin{proof}
  Since $u$ is a unit, $x \mapsto ux$ is a bijection of $\mathbb{Z}_n$, so $\varphi_u$ is a bijection of $V$.
  Every bijection of the vertex set of a complete graph is an automorphism of it.

  For a pair $\{x, y\} \in S$ we have $\varphi_u(\{x + i, y + i\}) = \{ux + ui, uy + ui\} = u\{x, y\} + ui$, so $\varphi_u(S + i) = uS + ui = S + ui$.
  For the remaining edge $\varphi_u(\{i, \infty\}) = \{ui, \infty\}$.
  Thus $\varphi_u(F_i) = (S + ui) \cup \{\{ui, \infty\}\} = F_{ui}$.
  In particular $\varphi_u(F_0) = F_0$ because $u \cdot 0 = 0$, so $\varphi_u(F_0 \cup F_d) = F_0 \cup F_{ud}$, and since $\varphi_u$ is an automorphism of $K_{n + 1}$, the union $F_0 \cup F_d$ is a Hamilton cycle if and only if $F_0 \cup F_{ud}$ is.
\end{proof}

The edge at $\infty$ is where the choice of $0$ as the element not covered by $S$ matters.
If $S$ covered $\mathbb{Z}_n \setminus \{m\}$ instead, then $F_i$ would contain $\{m + i, \infty\}$, which $\varphi_u$ maps to $\{um + ui, \infty\}$, and this edge lies in $F_{ui}$ only when $um = m$.

For a subgroup $U$ of $\mathbb{Z}_n^{*}$ we write $G = \langle -1, U \rangle$ for the subgroup of $\mathbb{Z}_n^{*}$ generated by $-1$ and $U$, which acts on $\mathbb{Z}_n \setminus \{0\}$ by multiplication.

\begin{lemma}\label{lem:orbit-conditions}
  If $U$ is a subgroup of $\mathbb{Z}_n^{*}$, $S$ is a $U$-invariant starter in $\mathbb{Z}_n$ and $D \subseteq \mathbb{Z}_n \setminus \{0\}$ contains at least one element of every $G$-orbit, then the 1-factorization induced by $S$ is perfect if and only if $F_0 \cup F_d$ is a Hamilton cycle for every $d \in D$.
\end{lemma}
\begin{proof}
  Let $H = \{d \in \mathbb{Z}_n \setminus \{0\} \mid F_0 \cup F_d \; \text{is a Hamilton cycle}\}$.
  Since $F_0 \cup F_{-d} \cong F_0 \cup F_d$, $d \in H$ implies $-d \in H$.
  By \cref{lem:multiplier-automorphism}, $d \in H$ implies $ud \in H$ for every $u \in U$.
  Every element of $G$ is a product of $-1$ and elements of $U$, so $H$ is closed under multiplication by $G$ and is therefore a union of $G$-orbits.
  If $D \subseteq H$, then $H$ meets every $G$-orbit and hence $H = \mathbb{Z}_n \setminus \{0\}$.
  Since $F_i \cup F_j \cong F_0 \cup F_{j - i}$, this is exactly perfection.
  Conversely, if the 1-factorization is perfect, then $H = \mathbb{Z}_n \setminus \{0\}$ contains $D$.
\end{proof}

\begin{example}\label{ex:conditions-91}
  Continuing \cref{ex:orbits-91}, let $G = \langle -1, 16 \rangle$.
  The $G$-orbit of $1$ is
  \[ \{\pm 1, \pm 16, \pm 74\} = \{1, 16, 17, 74, 75, 90\}, \]
  so if $F_0 \cup F_1$ is a Hamilton cycle, then so is $F_0 \cup F_d$ for every
  \[ d \in \{1, 16, 17, 74, 75, 90\}. \]
\end{example}

Only the elements of $U$ come from automorphisms of the 1-factorization; the element $-1$ comes from translation, through the isomorphism $F_0 \cup F_d \cong F_0 \cup F_{-d}$.
Multiplication by $-1$ itself is an automorphism only when $-S = S$, and the next lemma shows that this happens only for the \emph{patterned starter} $P = \{\{x, -x\} \mid x \in \mathbb{Z}_n \setminus \{0\}\}$.

\begin{lemma}\label{lem:minus-one-patterned}
  If $S$ is a starter in $\mathbb{Z}_n$ with $-S = S$, then $S = P$.
\end{lemma}
\begin{proof}
  Suppose $\{x, y\} \in S$ with $y \ne -x$.
  Then $\{-x, -y\} = -\{x, y\}$ is also in $S$.
  It is a different pair: $-x \ne x$ because $n$ is odd and $x \ne 0$, and $-x \ne y$ by assumption.
  But both pairs have the differences $\pm(x - y)$, which contradicts the fact that every difference occurs exactly once in $S$.
  Thus every pair of $S$ has the form $\{x, -x\}$, and since $S$ has $\frac{n - 1}{2}$ pairs, $S = P$.
\end{proof}

\subsection{Admissible Subgroups}\label{sec:admissible}

Not every subgroup $U$ can lead to a perfect 1-factorization.
To see which ones can, we need a general fact about symmetries of perfect 1-factorizations.

By \cref{lem:multiplier-automorphism}, $\varphi_u$ is an automorphism of the 1-factorization induced by a $U$-invariant starter for every $u \in U$.

\begin{lemma}[{\cite[Remark~1.3(d) on p.~311]{ihrig:StructureSymmetryGroups:1989}}]\label{lem:two-fixed-points}
  If $\alpha$ is a non-identity automorphism of a perfect 1-factorization, then $\alpha$ fixes at most two vertices.
\end{lemma}

\begin{lemma}\label{lem:fixed-point-free}
  If $U$ is a subgroup of $\mathbb{Z}_n^{*}$ and $S$ is a $U$-invariant starter in $\mathbb{Z}_n$ that induces a perfect 1-factorization, then $\gcd(u - 1, n) = 1$ for every $u \in U \setminus \{1\}$.
\end{lemma}
\begin{proof}
  Let $u \in U \setminus \{1\}$.
  Since $\gcd(u - 1, n) > 1$ holds exactly when some prime $p$ divides both $n$ and $u - 1$, it is equivalent to show that $u \not\equiv 1 \pmod{p}$ for every prime $p$ dividing $n$.
  The automorphism $\varphi_u$ is not the identity, since $\varphi_u(1) = u \ne 1$, and it fixes the vertices $0$ and $\infty$.
  By \cref{lem:two-fixed-points} $\varphi_u$ fixes no other vertex, so $ux = x$, that is $(u - 1)x \equiv 0 \pmod{n}$, has no solution when $x \ne 0$.
  If $u \equiv 1 \pmod{p}$ for some prime $p$ dividing $n$, then $x = \frac{n}{p}$ is such a solution, because $\frac{n \left( u - 1 \right)}{p} \equiv 0 \pmod{n}$.
  Thus $u \not\equiv 1 \pmod{p}$ for every prime $p$ dividing $n$.
\end{proof}

\begin{example}\label{ex:residues-91}
  In \cref{ex:orbits-91} we have $n = 91 = 7 \cdot 13$, and the condition of \cref{lem:fixed-point-free} holds for both non-identity elements of $U = \langle 16 \rangle$: $16 \equiv 2 \pmod{7}$ and $16 \equiv 3 \pmod{13}$, and $74 \equiv 4 \pmod{7}$ and $74 \equiv 9 \pmod{13}$.
\end{example}

The element $-1$ satisfies the condition of \cref{lem:fixed-point-free} for every odd $n$, since $-1 - 1 = -2$ and $\gcd(2, n) = 1$.
If $-1 \in U$, then $S$ is invariant under $-1$, so $S = P$ by \cref{lem:minus-one-patterned}.
So to rule out $-1 \in U$ for composite $n$, it suffices that $P$ does not induce a perfect 1-factorization, which we show below.
Without this argument $-1$ would look like a free additional symmetry: it leaves $G$ and hence the number of conditions unchanged, but it halves the number of orbits of pairs that form a starter and so shrinks the search space further.

\begin{lemma}\label{lem:patterned-not-perfect}
  If $n$ is composite, then the patterned starter $P$ does not induce a perfect 1-factorization of $K_{n + 1}$.
\end{lemma}
\begin{proof}
  Let $F_i$ be the factors induced by $P$, let $d \in \mathbb{Z}_n \setminus \{0\}$ and let $A = \{0, d, 2d, \ldots\}$ be the additive subgroup of $\mathbb{Z}_n$ generated by $d$.
  In $F_0$ the neighbour of a vertex $y \ne 0$ is $-y$, and the neighbour of $0$ is $\infty$.
  In $F_d = (P + d) \cup \{\{d, \infty\}\}$ the neighbour of a vertex $y \ne d$ is $2d - y$, and the neighbour of $d$ is $\infty$.
  So every neighbour in $F_0 \cup F_d$ of a vertex of $A \cup \{\infty\}$ lies again in $A \cup \{\infty\}$.
  Hence the cycle of $F_0 \cup F_d$ through $\infty$ uses only vertices of $A \cup \{\infty\}$.

  Since $n$ is composite, it has a prime divisor $p < n$.
  For $d = p$ the subgroup $A$ has $\frac{n}{p} < n$ elements, so the cycle through $\infty$ has at most $\frac{n}{p} + 1 < n + 1$ vertices and $F_0 \cup F_p$ is not a Hamilton cycle.
\end{proof}

\begin{lemma}\label{lem:minus-one}
  If $n$ is composite, $U$ is a subgroup of $\mathbb{Z}_n^{*}$ and $S$ is a $U$-invariant starter in $\mathbb{Z}_n$ that induces a perfect 1-factorization, then $-1 \notin U$.
\end{lemma}
\begin{proof}
  If $-1 \in U$, then $-S = S$, so $S = P$ by \cref{lem:minus-one-patterned}, which contradicts \cref{lem:patterned-not-perfect}.
\end{proof}

\Cref{lem:fixed-point-free,lem:minus-one} give two necessary conditions on $U$.
We say that a subgroup $U$ of $\mathbb{Z}_n^{*}$ is \emph{admissible} if $\gcd(u - 1, n) = 1$ for every $u \in U \setminus \{1\}$ and $-1 \notin U$.
So when $n$ is composite, every subgroup $U$ for which a $U$-invariant starter induces a perfect 1-factorization is admissible.
The next lemma describes the admissible subgroups in terms that are easy to check.

\begin{lemma}\label{lem:admissible}
  If $n \ge 3$ is odd and $U$ is a subgroup of $\mathbb{Z}_n^{*}$, then:
  \begin{enumerate}[label = (\alph*)]
    \item $U$ is admissible if and only if $|U|$ is odd and, for every prime $p$ dividing $n$, no element of $U \setminus \{1\}$ is congruent to $1$ modulo $p$;
    \item if $U$ is admissible, then $U$ is cyclic and $|U|$ divides $p - 1$ for every prime $p$ dividing $n$.
  \end{enumerate}
\end{lemma}
\begin{proof}\hspace{0pt}
  \begin{enumerate}[label = (\alph*)]
    \item
          As noted in the proof of \cref{lem:fixed-point-free}, the condition $\gcd(u - 1, n) = 1$ for all $u \in U \setminus \{1\}$ is the same as the condition on primes in the statement.
          So it remains to show that, given this condition, $-1 \notin U$ holds if and only if $|U|$ is odd.

          If $|U|$ is odd, then $U$ contains no element of order $2$, and $-1$ has order $2$ since $n \ge 3$, thus $-1 \notin U$.
          Conversely, suppose $-1 \notin U$ but $|U|$ is even.
          By Cauchy's theorem $U$ contains an element $u$ of order $2$, and $u \ne -1$.
          Let $p^a$ be a prime power exactly dividing $n$, that is, $p^a$ divides $n$ but $p^{a + 1}$ does not.
          Since $u^2 = 1$, the prime power $p^a$ divides $(u - 1)(u + 1)$.
          By assumption $u \not\equiv 1 \pmod{p}$, so $p$ does not divide $u - 1$, and therefore $p^a$ divides $u + 1$.
          As this holds for every prime power exactly dividing $n$, the Chinese remainder theorem gives $u \equiv -1 \pmod{n}$, a contradiction.
    \item
          Fix a prime $p$ dividing $n$ and consider reduction modulo $p$, a homomorphism from $U$ to $\mathbb{Z}_p^{*}$.
          Its kernel is $\{u \in U \mid u \equiv 1 \pmod{p}\} = \{1\}$ by (a), so $U$ is isomorphic to a subgroup of $\mathbb{Z}_p^{*}$.
          The group $\mathbb{Z}_p^{*}$ is cyclic of order $p - 1$, so $U$ is cyclic and $|U|$ divides $p - 1$.\qedhere
  \end{enumerate}
\end{proof}

In particular, an admissible $U \ne \{1\}$ exists only if the numbers $p - 1$ for the primes $p$ dividing $n$ have a common odd divisor greater than $1$.
The two conditions of the admissible subgroups also follow from a theorem of Ihrig on perfect 1-factorizations induced by starters~\cite[Theorem~4.1(b) on p.~321]{ihrig:StructureSymmetryGroups:1989}.
For such a 1-factorization of $K_{n + 1}$ that is not induced by the patterned starter and in which every vertex is fixed by a non-identity automorphism, as is the case when $S$ is $U$-invariant for some $U \ne \{1\}$, it shows that the automorphisms fixing $0$ and $\infty$ are the multiplications by the elements of a subgroup of $\mathbb{Z}_n^{*}$ of odd order whose non-identity elements fix only $0$.
\Cref{lem:fixed-point-free,lem:minus-one} prove these conditions for $U$ directly.

\subsection{Orbits}\label{sec:orbits}

We now show that for an admissible subgroup all the orbits involved have the same size, so that the reduction of \cref{sec:multiplier} can be counted exactly.
This follows from the fact that the actions involved are free, where a group $H$ acts \emph{freely} on a set $X$ if $hx = x$ for some $x \in X$ implies $h = 1$.

\begin{lemma}\label{lem:free}
  If $U$ is an admissible subgroup of $\mathbb{Z}_n^{*}$, then:
  \begin{enumerate}[label = (\alph*)]
    \item $G = U \cup (-U)$ with $U \cap (-U) = \emptyset$, so $|G| = 2|U|$;
    \item $U$ and $G$ act freely on $\mathbb{Z}_n \setminus \{0\}$;
    \item $U$ acts freely on the set of pairs of elements of $\mathbb{Z}_n \setminus \{0\}$;
    \item every $U$-orbit in $\mathbb{Z}_n \setminus \{0\}$ has $|U|$ elements, every $G$-orbit has $2|U|$ elements, and every $U$-orbit of pairs has $|U|$ pairs.
  \end{enumerate}
\end{lemma}
\begin{proof}\hspace{0pt}
  \begin{enumerate}[label = (\alph*)]
    \item
          The set $U \cup (-U)$ contains $-1$ and $U$, and it is closed under multiplication because $\mathbb{Z}_n^{*}$ is commutative and ${(-1)}^2 = 1$.
          So it is the subgroup generated by $-1$ and $U$.
          If $u = -v$ for some $u, v \in U$, then $-1 = u v^{-1} \in U$, which is impossible, so $U$ and $-U$ are disjoint.
    \item
          Let $d \in \mathbb{Z}_n \setminus \{0\}$.
          If $ud = d$ for some $u \in U \setminus \{1\}$, then $(u - 1)d \equiv 0 \pmod{n}$ with $\gcd(u - 1, n) = 1$, so $d = 0$, a contradiction.
          If $-ud = d$ for some $u \in U$, then ${(-u)}^{|U|} = -u^{|U|} = -1$ also fixes $d$, because $|U|$ is odd by \cref{lem:admissible}~(a).
          So $2d = 0$, and since $n$ is odd, $d = 0$, again a contradiction.
    \item
          Suppose $u\{x, y\} = \{x, y\}$ for some $u \in U \setminus \{1\}$.
          By (b), $ux \ne x$, so $ux = y$ and $uy = x$.
          Then $u^2 x = x$, so $u^2 = 1$ by (b), and $u$ has order $2$.
          This is impossible because $|U|$ is odd.
    \item
          By the orbit-stabilizer theorem, if a finite group $H$ acts on a set $X$, then $|Hx| \cdot |H_x| = |H|$ for every $x \in X$, where $H_x = \{h \in H \mid hx = x\}$ is the stabilizer of $x$.
          Let $d \in \mathbb{Z}_n \setminus \{0\}$.
          By (b) the stabilizers $U_d$ and $G_d$ are trivial, so $|Ud| = |U|$, and $|Gd| = |G| = 2|U|$ by (a).
          By (c) the stabilizer in $U$ of a pair of elements of $\mathbb{Z}_n \setminus \{0\}$ is also trivial, so its $U$-orbit has $|U|$ pairs.\qedhere
  \end{enumerate}
\end{proof}

Part (c) is what lets us treat a $U$-invariant starter as a collection of whole orbits of pairs, each containing exactly $|U|$ pairs.
Not every orbit of pairs can be part of a starter, because the pairs of a starter are disjoint.
For example, if $|U| = 3$ and $U = \{1, u, u^2\}$, the orbit of $\{x, ux\}$ consists of $\{x, ux\}$, $\{ux, u^2x\}$ and $\{u^2x, x\}$, which pairwise share an element.
The next lemma shows that this is the only obstruction.

\begin{lemma}\label{lem:orbit-disjoint}
  If $U$ is an admissible subgroup of $\mathbb{Z}_n^{*}$ and $x, y \in \mathbb{Z}_n \setminus \{0\}$ are distinct, then:
  \begin{enumerate}[label = (\alph*)]
    \item the pairs in the $U$-orbit of $\{x, y\}$ are pairwise disjoint if and only if $y \notin Ux$;
    \item if $y \notin Ux$, then the orbit covers the $2|U|$ elements of $Ux \cup Uy$, and the differences $\pm(ux - uy)$, $u \in U$, of its pairs are exactly the $2|U|$ elements of the $G$-orbit of $x - y$, each occurring once.
  \end{enumerate}
\end{lemma}
\begin{proof}\hspace{0pt}
  \begin{enumerate}[label = (\alph*)]
    \item
          If $y = vx$ for some $v \in U$, then $v \ne 1$ because $y \ne x$, so $\{x, y\}$ and $v\{x, y\} = \{vx, v^2x\}$ are different pairs by \cref{lem:free}~(c), and they share the element $vx = y$.

          Conversely, let $y \notin Ux$ and let $u \ne v$ be elements of $U$.
          We have $ux \ne vx$ and $uy \ne vy$ by \cref{lem:free}~(b).
          If $ux = vy$, then $y = v^{-1}u x \in Ux$, and similarly $uy = vx$ is impossible.
          So $\{ux, uy\}$ and $\{vx, vy\}$ are disjoint.
    \item
          Since $y \notin Ux$, the orbits $Ux$ and $Uy$ are disjoint, and each has $|U|$ elements by \cref{lem:free}~(d).
          The pair $u\{x, y\}$ has the differences $\pm u(x - y)$.
          As $u$ runs over $U$ these are the elements $g(x - y)$ for $g \in U \cup (-U) = G$.
          Since $x - y \ne 0$ and $G$ acts freely, these $2|U|$ elements are distinct.\qedhere
  \end{enumerate}
\end{proof}

We can now count everything exactly.

\begin{lemma}\label{lem:counts}
  If $U$ is an admissible subgroup of $\mathbb{Z}_n^{*}$, then:
  \begin{enumerate}[label = (\alph*)]
    \item the group $G$ has exactly $\frac{n - 1}{2|U|}$ orbits on $\mathbb{Z}_n \setminus \{0\}$;
    \item a set $S$ of pairs is a $U$-invariant starter if and only if it is the union of $\frac{n - 1}{2|U|}$ orbits of pairs $\{x, y\}$ with $y \notin Ux$, such that every element of $\mathbb{Z}_n \setminus \{0\}$ is covered by a pair of exactly one of the chosen orbits and every $G$-orbit of $\mathbb{Z}_n \setminus \{0\}$ is the set of differences of exactly one of the chosen orbits.
  \end{enumerate}
\end{lemma}
\begin{proof}\hspace{0pt}
  \begin{enumerate}[label = (\alph*)]
    \item
          The set $\mathbb{Z}_n \setminus \{0\}$ has $n - 1$ elements, and by \cref{lem:free}~(d) it splits into $G$-orbits of $2|U|$ elements each.
    \item
          Let $S$ be a $U$-invariant starter.
          It is a union of orbits of pairs, and since its pairs are disjoint, by \cref{lem:orbit-disjoint}~(a) each of these orbits is the orbit of a pair $\{x, y\}$ with $y \notin Ux$.
          Its pairs cover every element of $\mathbb{Z}_n \setminus \{0\}$ once, so every element is covered by a pair of exactly one of the orbits.
          Its differences cover every element of $\mathbb{Z}_n \setminus \{0\}$ once, and by \cref{lem:orbit-disjoint}~(b) each orbit contributes exactly one whole $G$-orbit of differences, so every $G$-orbit comes from exactly one of the orbits.
          The number of orbits is the number of pairs, $\frac{n - 1}{2}$, divided by the orbit size $|U|$, that is $\frac{n - 1}{2|U|}$.

          Conversely, a union of orbits with these properties is $U$-invariant, its pairs are disjoint and cover $\mathbb{Z}_n \setminus \{0\}$, and by \cref{lem:orbit-disjoint}~(b) its differences cover every element of $\mathbb{Z}_n \setminus \{0\}$ exactly once.
          So it is a starter.\qedhere
  \end{enumerate}
\end{proof}

For starters in $\mathbb{F}_q$, the conditions of \cref{lem:orbit-disjoint}~(a) and \cref{lem:counts}~(b) appear in~\cite[Lemma~2.1 on p.~407]{dinitz:NewPerfectOnefactorizations:1989}, where they are stated without proof.
By part (a), \cref{lem:orbit-conditions} leaves exactly $\frac{n - 1}{2|U|}$ Hamilton cycle conditions for a $U$-invariant starter.
Part (b) turns the search for a $U$-invariant starter into an \emph{exact cover} problem.
The items to be covered are the $n - 1$ elements of $\mathbb{Z}_n \setminus \{0\}$ together with the $\frac{n - 1}{2|U|}$ orbits of $G$, and each orbit of pairs $\{x, y\}$ with $y \notin Ux$ covers the $2|U|$ elements of $Ux \cup Uy$ and one $G$-orbit.
A $U$-invariant starter is exactly a choice of orbits of pairs that covers every item once.

\subsection{Isomorphism}\label{sec:isomorphism}

Multiplier symmetry also gives a simple way to tell 1-factorizations apart, using the following result of Ihrig.

\begin{lemma}[{\cite[Corollary~4.4 on p.~141]{ihrig:ConditionsEnsuringStarter:2003}}]\label{lem:ihrig-isomorphism}
  If $n \ne 27$, $S$ and $S'$ are starters in $\mathbb{Z}_n$ and the union of any three distinct factors of the 1-factorization induced by $S$ is connected, then the 1-factorizations induced by $S$ and $S'$ are isomorphic if and only if $S' = cS$ for some $c \in \mathbb{Z}_n^{*}$.
\end{lemma}

\begin{lemma}\label{lem:non-isomorphic}
  If $n \ne 27$, $a, b \in \mathbb{Z}_n^{*}$ satisfy $ab \ne 1$ and $\gcd(ab - 1, n) > 1$, and $S_a$ and $S_b$ are starters in $\mathbb{Z}_n$ with $aS_a = S_a$ and $bS_b = S_b$ that induce perfect 1-factorizations, then these 1-factorizations are not isomorphic.
\end{lemma}
\begin{proof}
  The union of any two factors of a perfect 1-factorization is a Hamilton cycle, so the union of any three distinct factors is connected.
  Suppose that the two 1-factorizations are isomorphic.
  By \cref{lem:ihrig-isomorphism}, $S_a = cS_b$ for some $c \in \mathbb{Z}_n^{*}$.
  Then $bS_a = bcS_b = c(bS_b) = cS_b = S_a$, so $S_a$ is invariant under the subgroup $\langle a, b \rangle$ of $\mathbb{Z}_n^{*}$.
  This subgroup contains $ab \ne 1$ with $\gcd(ab - 1, n) > 1$, which contradicts \cref{lem:fixed-point-free}.
\end{proof}

\subsection{The Case \texorpdfstring{$n = 91$}{n = 91}}\label{sec:91}

We now determine all admissible subgroups for $n = 91$.
By the Chinese remainder theorem, the map $u \mapsto (u \bmod 7, u \bmod 13)$ identifies $\mathbb{Z}_{91}^{*}$ with $\mathbb{Z}_7^{*} \times \mathbb{Z}_{13}^{*}$, a group of order $6 \cdot 12 = 72$.
By \cref{lem:admissible}, the order of an admissible subgroup is odd and divides both $6$ and $12$, so it is $1$ or $3$.

The elements of order $3$ in $\mathbb{Z}_7^{*}$ are $2$ and $4$, and those in $\mathbb{Z}_{13}^{*}$ are $3$ and $9$.
So the elements $u \ne 1$ with $u^3 = 1$ in $\mathbb{Z}_{91}^{*}$ are the eight pairs $(a, b) \ne (1, 1)$ with $a \in \{1, 2, 4\}$ and $b \in \{1, 3, 9\}$, and they form four subgroups of order $3$:
\begin{table}[H]
  \centering
  \begin{tabular}{lll}
    \toprule
    subgroup        & residues modulo $7$ and $13$ & admissible                  \\
    \midrule
    $\{1, 16, 74\}$ & $(1, 1), (2, 3), (4, 9)$     & yes                         \\
    $\{1, 9, 81\}$  & $(1, 1), (2, 9), (4, 3)$     & yes                         \\
    $\{1, 22, 29\}$ & $(1, 1), (1, 9), (1, 3)$     & no, $22 \equiv 1 \pmod{7}$  \\
    $\{1, 53, 79\}$ & $(1, 1), (4, 1), (2, 1)$     & no, $53 \equiv 1 \pmod{13}$ \\
    \bottomrule
  \end{tabular}
\end{table}
A subgroup is admissible exactly when its non-identity elements are $\ne 1$ in both coordinates.
So there are exactly two admissible subgroups of order $3$, $\langle 16 \rangle$ and $\langle 9 \rangle$, which correspond to the two ways of matching the elements of order $3$ modulo $7$ with those modulo $13$.

\begin{example}[Inadmissible Subgroup]\label{ex:22}
  To see concretely what goes wrong for the other two subgroups, take $U = \langle 22 \rangle$.
  Since $22 - 1 = 21 = 3 \cdot 7$, multiplication by $22$ fixes every multiple of $13$, for example $22 \cdot 13 = 286 \equiv 13 \pmod{91}$.
  So $\varphi_{22}$ fixes the three vertices $0$, $13$ and $\infty$, and by \cref{lem:two-fixed-points} no $\langle 22 \rangle$-invariant starter induces a perfect 1-factorization.
  The orbit counts of \cref{lem:counts} also fail: under $\langle -1, 22 \rangle$ the six multiples of $13$ form three orbits $\{\pm 13\}$, $\{\pm 26\}$, $\{\pm 39\}$ of size $2$, and the remaining $84$ elements form $14$ orbits of size $6$, giving $17$ orbits instead of $15$.
\end{example}

\begin{proof}[Proof of \cref{thm:main}]
  Let $U = \langle 16 \rangle$, which is admissible, and let $G = \langle -1, 16 \rangle$.
  One set of representatives of the $15$ orbits of $G$ is
  \[
    D_{16} = \{1, 2, 3, 4, 5, 7, 8, 9, 10, 13, 14, 15, 19, 20, 25\}.
  \]
  Let $S$ be the union of the $U$-orbits of the $15$ pairs
  \begin{gather*}
    \{1, 41\}, \{2, 40\}, \{4, 25\}, \{5, 79\}, \{7, 26\}, \{8, 35\}, \{9, 17\}, \{10, 45\}, \\
    \{11, 44\}, \{15, 62\}, \{20, 59\}, \{27, 39\}, \{28, 60\}, \{33, 88\}, \{49, 55\}.
  \end{gather*}
  For each of these pairs $\{x, y\}$ we have $y \notin Ux$, so by \cref{lem:orbit-disjoint} its orbit consists of three disjoint pairs that cover $Ux \cup Uy$, and the differences of these pairs form the $G$-orbit of $x - y$.
  The $30$ orbits $Ux$ and $Uy$ obtained in this way are distinct, so they partition $\mathbb{Z}_{91} \setminus \{0\}$.
  The differences $x - y$ of the $15$ pairs lie in distinct $G$-orbits, so every $G$-orbit is the set of differences of exactly one of the chosen orbits.
  By \cref{lem:counts}~(b), $S$ is a $U$-invariant starter; it is the set $S_{16}$ listed in \cref{app:starters}.
  A direct computation shows that for each of the $15$ values $d \in D_{16}$ the cycle of $F_0 \cup F_d$ through $\infty$ contains all $92$ vertices.
  By \cref{lem:orbit-conditions}, the 1-factorization of $K_{92}$ induced by $S$ is perfect.

  The other admissible subgroup $\langle 9 \rangle$ also works.
  The same argument, with the set of representatives
  \[
    D_9 = \{1, 2, 3, 4, 5, 6, 7, 8, 12, 13, 14, 15, 16, 23, 24\}
  \]
  of the orbits of $\langle -1, 9 \rangle$ in place of $D_{16}$, shows that the union of the $\langle 9 \rangle$-orbits of the pairs
  \begin{gather*}
    \{1, 2\}, \{3, 5\}, \{4, 8\}, \{6, 12\}, \{7, 15\}, \{10, 43\}, \{13, 83\}, \{14, 66\}, \\
    \{16, 65\}, \{20, 79\}, \{24, 49\}, \{28, 55\}, \{30, 75\}, \{37, 59\}, \{46, 58\}
  \end{gather*}
  is a $\langle 9 \rangle$-invariant starter, the set $S_9$ listed in \cref{app:starters}, that induces a perfect 1-factorization of $K_{92}$.

  Finally, $9 \cdot 16 \equiv 53 \pmod{91}$ and $\gcd(53 - 1, 91) = 13$, so the two 1-factorizations are not isomorphic by \cref{lem:non-isomorphic}.
\end{proof}

We found $S_{16}$ and $S_9$ by a backtracking search over the exact cover problem of \cref{sec:orbits} that abandons a partial solution as soon as its pairs close a cycle missing a vertex in some $F_0 \cup F_d$ with $d \in D_{16}$ or $d \in D_9$ respectively.

\section{Future Work}\label{sec:future}

By \cref{lem:fixed-point-free,lem:minus-one}, the method applies to $K_{n + 1}$ only if $\mathbb{Z}_n^{*}$ has an admissible subgroup $U \ne \{1\}$.
Among the other eight orders up to $100$ for which no perfect 1-factorization is known, none has one by \cref{lem:admissible}.
For example, the method says nothing about $K_{64}$, since $3$ divides $63$ and every element of order $3$ in $\mathbb{Z}_{63}^{*}$ is congruent to $1$ modulo $3$.

The table below lists the orders $2t \le 1000$ for which no perfect 1-factorization is known~\cite{pike:PerfectOnefactorisation:2019,rosa:Perfect1factorizations:2019,wanless:Perfect1factorisations:} and for which $\mathbb{Z}_{2t - 1}^{*}$ has an admissible subgroup $U \ne \{1\}$, together with the largest possible $|U|$ and the resulting number of conditions $c = \frac{2t - 2}{2|U|}$ from \cref{lem:counts}.
The orders $50$, $170$, $344$ and $530$ also meet the criterion but are not listed because perfect 1-factorizations of these orders are already known~\cite{rosa:Perfect1factorizations:2019,wanless:Perfect1factorisations:}; for example, $529 = 23^2$ admits $|U| = 11$ in $\mathbb{Z}_{529}$ and $|U| = 33$ in $\mathbb{F}_{529}$, the latter with $c = 8$.
\begin{table}[H]
  \centering
  \begin{tabular}[t]{rlrr}
    \toprule
    $2t$  & $2t - 1$       & $|U|$ & $c$   \\
    \midrule
    $92$  & $7 \cdot 13$   & $3$   & $15$  \\
    $248$ & $13 \cdot 19$  & $3$   & $41$  \\
    $260$ & $7 \cdot 37$   & $3$   & $43$  \\
    $342$ & $11 \cdot 31$  & $5$   & $34$  \\
    $404$ & $13 \cdot 31$  & $3$   & $67$  \\
    $428$ & $7 \cdot 61$   & $3$   & $71$  \\
    $452$ & $11 \cdot 41$  & $5$   & $45$  \\
    $470$ & $7 \cdot 67$   & $3$   & $78$  \\
    $512$ & $7 \cdot 73$   & $3$   & $85$  \\
    $560$ & $13 \cdot 43$  & $3$   & $93$  \\
    $590$ & $19 \cdot 31$  & $3$   & $98$  \\
    $638$ & $7^2 \cdot 13$ & $3$   & $106$ \\
    \bottomrule
  \end{tabular}
  \qquad
  \begin{tabular}[t]{rlrr}
    \toprule
    $2t$  & $2t - 1$       & $|U|$ & $c$   \\
    \midrule
    $672$ & $11 \cdot 61$  & $5$   & $67$  \\
    $680$ & $7 \cdot 97$   & $3$   & $113$ \\
    $704$ & $19 \cdot 37$  & $9$   & $39$  \\
    $722$ & $7 \cdot 103$  & $3$   & $120$ \\
    $764$ & $7 \cdot 109$  & $3$   & $127$ \\
    $782$ & $11 \cdot 71$  & $5$   & $78$  \\
    $872$ & $13 \cdot 67$  & $3$   & $145$ \\
    $890$ & $7 \cdot 127$  & $3$   & $148$ \\
    $932$ & $7^2 \cdot 19$ & $3$   & $155$ \\
    $950$ & $13 \cdot 73$  & $3$   & $158$ \\
    $962$ & $31^2$         & $15$  & $32$  \\
    \bottomrule
  \end{tabular}
\end{table}
The choice of $U$ enters in two opposite ways.
A larger $U$ gives fewer conditions, which makes each candidate starter more likely to work, but it also leaves fewer $U$-invariant starters, and for large $|U|$ there may be none that works.
A smaller $U$ has more candidates but more conditions.
For example, $K_{962}$ admits $|U| \in \{3, 5, 15\}$ with $c = 160$, $96$ and $32$, and $K_{704}$ admits $|U| \in \{3, 9\}$ with $c = 117$ and $39$.
Understanding this trade-off, and the actual number of $U$-invariant starters as a function of $|U|$, would help to decide which orders are worth attacking next.

\appto{\bibsetup}{\emergencystretch=1em}
\printbibliography%

@book{colbourn:HandbookCombinatorialDesigns:2006,
  title = {Handbook of Combinatorial Designs},
  editor = {Colbourn, Charles J. and Dinitz, Jeffrey H.},
  date = {2006-11-02},
  edition = {2},
  publisher = {{Chapman and Hall/CRC}},
  doi = {10.1201/9781420010541},
  url = {https://www.taylorfrancis.com/books/9781420010541},
  isbn = {978-0-429-13848-5},
  langid = {english}
}

@article{dinitz:NewPerfectOnefactorizations:1989,
  title = {Some new perfect one-factorizations from starters in finite fields},
  author = {Dinitz, Jeffrey H. and Stinson, Douglas R.},
  date = {1989},
  journaltitle = {Journal of Graph Theory},
  volume = {13},
  number = {4},
  pages = {405--415},
  issn = {1097-0118},
  doi = {10.1002/jgt.3190130404},
  url = {https://onlinelibrary.wiley.com/doi/abs/10.1002/jgt.3190130404},
  langid = {english}
}

@article{ihrig:ConditionsEnsuringStarter:2003,
  title = {Conditions ensuring starter induced 1-factorizations are not isomorphic},
  author = {Ihrig, Edwin C.},
  date = {2003-01},
  journaltitle = {Journal of Combinatorial Designs},
  volume = {11},
  number = {2},
  pages = {124--143},
  publisher = {Wiley},
  issn = {1063-8539},
  doi = {10.1002/jcd.10023},
  url = {https://doi.org/10.1002/jcd.10023},
}

@article{ihrig:PerfectOnefactorizationK50:1987,
  title = {A perfect one-factorization of $K_{50}$},
  author = {Ihrig, Edwin C. and Seah, Eric and Stinson, Douglas R.},
  date = {1987},
  journaltitle = {Journal of Combinatorial Mathematics and Combinatorial Computing},
  shortjournal = {J. Combin. Math. Combin. Comput.},
  volume = {1},
  pages = {217--219},
  issn = {0835-3026},
  url = {https://combinatorialpress.com/article/jcmcc/Volume%201/vol-001-paper%2017.pdf},
  langid = {english}
}

@article{ihrig:StructureSymmetryGroups:1989,
  title = {The structure of the symmetry groups of perfect 1-factorizations of $K_{2n}$},
  author = {Ihrig, Edwin C.},
  date = {1989},
  journaltitle = {Journal of Combinatorial Theory, Series B},
  volume = {47},
  number = {3},
  pages = {307--329},
  doi = {10.1016/0095-8956(89)90031-2}
}

@article{pike:PerfectOnefactorisation:2019,
  title = {A perfect one-factorisation of $K_{56}$},
  author = {Pike, David A.},
  date = {2019},
  journaltitle = {Journal of Combinatorial Designs},
  volume = {27},
  number = {6},
  pages = {386--390},
  issn = {1520-6610},
  doi = {10.1002/jcd.21653},
  url = {https://onlinelibrary.wiley.com/doi/abs/10.1002/jcd.21653},
  langid = {english}
}

@article{rosa:Perfect1factorizations:2019,
  title = {Perfect 1-factorizations},
  author = {Rosa, Alexander},
  date = {2019-06-01},
  journaltitle = {Mathematica Slovaca},
  volume = {69},
  number = {3},
  pages = {479--496},
  publisher = {De Gruyter},
  issn = {1337-2211},
  doi = {10.1515/ms-2017-0241},
  url = {https://www.degruyterbrill.com/document/doi/10.1515/ms-2017-0241/html?lang=en},
  langid = {english}
}

@online{wanless:Perfect1factorisations:,
  title = {Perfect 1-factorisations},
  author = {Wanless, Ian M.},
  urldate  = {2026-10-09},
  url = {https://users.monash.edu.au/~iwanless/data/P1F/newP1F.html}
}

@article{wolfe:PerfectOnefactorizationK52:2009,
  title = {A perfect one-factorization of $K_{52}$},
  author = {Wolfe, Adam J.},
  date = {2009},
  journaltitle = {Journal of Combinatorial Designs},
  volume = {17},
  number = {2},
  pages = {190--196},
  issn = {1520-6610},
  doi = {10.1002/jcd.20194},
  url = {https://onlinelibrary.wiley.com/doi/abs/10.1002/jcd.20194},
  langid = {english}
}

@incollection{kotzig:HamiltonGraphsHamilton:1964,
  title = {Hamilton graphs and Hamilton circuits},
  booktitle = {Theory of Graphs and Its Applications: Proceedings of the Symposium Held in Smolenice in June 1963},
  author = {Kotzig, Anton},
  editor = {Fiedler, Miroslav},
  date = {1964},
  pages = {63--82},
  publisher = {Publishing House of the Czechoslovak Academy of Sciences},
  location = {Prague}
}

\appendix
\crefalias{section}{appendix}

\section{Starters of \texorpdfstring{$K_{92}$}{K92}}\label{app:starters}

For reference we list the two starters of \cref{thm:main} in full, each as a set of $45$ pairs in $\mathbb{Z}_{91}$.
Every three consecutive pairs form one $U$-orbit $\{x, y\}$, $u\{x, y\}$, $u^2\{x, y\}$, where $u$ is the generator $16$ or $9$ of $U$.

\begin{align*}
  S_{16} = \{ & \{1, 41\}, \{16, 19\}, \{31, 74\}, \{2, 40\}, \{3, 32\}, \{48, 57\},    \\
              & \{4, 25\}, \{36, 64\}, \{23, 30\}, \{5, 79\}, \{80, 81\}, \{6, 22\},    \\
              & \{7, 26\}, \{21, 52\}, \{13, 63\}, \{8, 35\}, \{14, 37\}, \{42, 46\},   \\
              & \{9, 17\}, \{53, 90\}, \{29, 75\}, \{10, 45\}, \{69, 83\}, \{12, 54\},  \\
              & \{11, 44\}, \{67, 85\}, \{71, 86\}, \{15, 62\}, \{58, 82\}, \{18, 38\}, \\
              & \{20, 59\}, \{34, 47\}, \{24, 89\}, \{27, 39\}, \{68, 78\}, \{65, 87\}, \\
              & \{28, 60\}, \{50, 84\}, \{70, 72\}, \{33, 88\}, \{43, 73\}, \{51, 76\}, \\
              & \{49, 55\}, \{56, 61\}, \{66, 77\} \}
\end{align*}
\begin{align*}
  S_9 = \{ & \{1, 2\}, \{9, 18\}, \{71, 81\}, \{3, 5\}, \{27, 45\}, \{41, 61\},      \\
           & \{4, 8\}, \{36, 72\}, \{11, 51\}, \{6, 12\}, \{17, 54\}, \{31, 62\},    \\
           & \{7, 15\}, \{44, 63\}, \{21, 32\}, \{10, 43\}, \{23, 90\}, \{25, 82\},  \\
           & \{13, 83\}, \{19, 26\}, \{52, 80\}, \{14, 66\}, \{35, 48\}, \{42, 68\}, \\
           & \{16, 65\}, \{39, 53\}, \{22, 78\}, \{20, 79\}, \{74, 89\}, \{29, 73\}, \\
           & \{24, 49\}, \{34, 77\}, \{33, 56\}, \{28, 55\}, \{40, 70\}, \{84, 87\}, \\
           & \{30, 75\}, \{38, 88\}, \{64, 69\}, \{37, 59\}, \{60, 76\}, \{47, 85\}, \\
           & \{46, 58\}, \{50, 67\}, \{57, 86\} \}
\end{align*}

For $S = S_{16}$ and $S = S_9$ the factors of the induced perfect 1-factorization of $K_{92}$ are $F_i = (S + i) \cup \{\{i, \infty\}\}$ for $i \in \mathbb{Z}_{91}$.

\end{document}